\documentclass[11pt,titlepage]{article}

\usepackage[latin1]{inputenc}
\usepackage{rotating, rotfloat}
\usepackage[round]{natbib} 
\usepackage{graphics, graphicx, hyperref, url, amsmath, subfigure, amsthm, amsfonts, enumerate, epstopdf, booktabs}
\hypersetup{
    pdfstartview={FitH},    % fits the width of the page to the window
    colorlinks=true,        % false: boxed links; true: colored links
    linkcolor=blue,         % color of internal links (change box color with linkbordercolor)
    citecolor=blue,         % color of links to bibliography
    filecolor=blue,         % color of file links
    urlcolor=blue           % color of external links
}
     \def\etal{{\itshape et al.}}
\def\beqn{\begin{eqnarray}} \def\eeqn{\end{eqnarray}}

\def\u{\textbf{\emph{u}}} \def\x{\textbf{\emph{x}}}  \def\z{\textbf{\emph{z}}}  \def\vv{\textbf{\emph{v}}}
\def\y{\textbf{\emph{y}}}

\newcommand{\bec}{\begin{center}}
\newcommand{\enc}{\end{center}}
\newcommand{\bee}{\begin{eqnarray*}}
\newcommand{\ene}{\end{eqnarray*}}
\newcommand{\beq}{\begin{equation}}
\newcommand{\eeq}{\end{equation}}

\begin{document}

\title{\bf Some results\\[1ex]
on the computing of Tukey's halfspace medain}
\vskip 5mm

\author {{Xiaohui Liu$^{a, b}$,   %\footnote{Corresponding author's email: csuliuxh912@gmail.com.}
        Shihua Luo$^{a, b}$,
        Yijun Zuo$^c$}\\ \\[1ex]
        {\em\footnotesize $^a$ School of Statistics, Jiangxi University of Finance and Economics, Nanchang, Jiangxi 330013, China}\\
         {\em\footnotesize $^b$ Research Center of Applied Statistics, Jiangxi University of Finance and Economics, Nanchang,}\\ {\em\footnotesize Jiangxi 330013, China}\\
         {\em\footnotesize $^c$ Department of Statistics and Probability, Michigan State University, East Lansing, MI, 48823, USA}\\
}

\maketitle

\begin{center}
{\sc Summary}
\end{center}

%\begin{center}
%\begin{minipage} {12cm}{ \bec {\sc Abstract  } \enc
%\vspace*{-2mm}
%{{\small \bf Abstract.}
Depth of the Tukey median is investigated for empirical distributions. A sharper upper bound is provided for this value for  data sets in general position. This bound is lower than the existing one in the literature, and more importantly derived under the \emph{fixed }sample size practical scenario. Several results obtained in this paper are interesting theoretically and useful as well to reduce the computational burden of the Tukey median practically when $p$ is large relative to large $n$.
%, such as $p = 5$ and $n = 5p$.
\vspace{2mm}

{\small {\bf\itshape Key words:} Data depth; Tukey median;  maximum halfspace depth; general position}
\vspace{2mm}

{\small {\bf2000 Mathematics Subject Classification Codes:} 62F10; 62F40; 62F35}
%\end{minipage}
%\end{center}

%-------------------------------------------------%
% Set base line skip
%-------------------------------------------------%
\setlength{\baselineskip}{1.5\baselineskip}

\vskip 0.1 in
\section{Introduction}
\paragraph{}
\vskip 0.1 in \label{Introduction}

For a univariate random sample $\mathcal{Z}^n = \{Z_1, Z_2, \cdots, Z_n\}$, let $Z_{(1)}, Z_{(2)}, \cdots, Z_{(n)}$ be the corresponding ordered statistics such that $Z_{(1)} \leq Z_{(2)} \leq \cdots \leq Z_{(n)}$. The sample median is then defined as
\begin{eqnarray*}
    \mu(\mathcal{Z}^n) = \frac{Z_{(\lfloor (n + 1) / 2 \rfloor)} + Z_{(\lfloor (n + 2) / 2 \rfloor)}}{2},
\end{eqnarray*}
where $\lfloor \cdot \rfloor$ denotes the floor function. In the literature, $\mu(\mathcal{Z}^n)$ is well known for its robustness. In fact, it has the highest possible breakdown point among its competitors and a bounded influence function. Hence, it usually serves as an alterative to the sample mean for data sets containing outliers (or with heavy tails).
\medskip

To extend the concept to multivariate setting is desirable.
The univariate definition %of the sample median $\mu(\mathcal{Z}^n)$
depends on the ordered statistics, whereas no natural ordering exists in the multivariate setting. It is thus nontrivial to extend it to multidimensional setting in $\mathcal{R}^p$ with $p \ge 2$. Several extensions, such as coordinatewise median \citep{Bic1964} and spatial median \citep{W09} and \citep{Bro1983}, etc., exist in the literature though. They lack the desirable affine equivariance, nevertheless. See \cite{Sma1990} for an early summary on multivariate medians.
\medskip

Observing the fact that $\mu(\mathcal{Z}^n)$ is the average of all points maximizing the function
\begin{eqnarray*}
    t \rightarrow \min\left\{P_n(Z \leq t), P_n(Z \ge t)\right\} = \min_{u \in \{ -1, 1\}} P_n(uZ \leq ut),
\end{eqnarray*}
where $Z$ is a univariate random variable and $P_n$ stands for the empirical distribution function,
\cite{Tuk1975} heuristically proposed a notion of depth function. He first defined the depth of a point $\x$ with respect to $\mathcal{X}^n$ as
\begin{eqnarray}
\label{DepthFunction}
    D(\x, \mathcal{X}^n) = \inf_{\u \in \{\z: \|\z\| = 1\}} P_n(\u^\top X \leq \u^\top \x),
\end{eqnarray}
where $\mathcal{X}^n = \{X_1, X_2, \cdots, X_n\}$ is a $p$-variate random sample, $X$ is a random variable in $R^p$ and $\|\cdot\|$ stands for Euclidean norm,
and then proposed to consider the average of such points that maximize $D(\x, \mathcal{X}^n)$ with respect to $\x$ as the multivariate median (more specially, hereafter Tukey median). That is,
\begin{eqnarray*}
    T^*(\mathcal{X}^n) = \textbf{Ave}\left\{\x \in \mathcal{R}^p : D(\x, \mathcal{X}^n) = \lambda^*(\mathcal{X}^n) \right\},
\end{eqnarray*}
where $\lambda^*(\mathcal{X}^n) = \sup_{\x} D(\x, \mathcal{X}^n)$. See also \cite{DG1992} (hereafter DG92) for details. But \emph{slightly differently}, DG92's discussions are based on an integer valued function, i.e., $nD(\x, \mathcal{X}^n)$, instead.
\medskip

%Since the above depth definition utilizes the closed halfspace $H_x:=\{y: u^\top y \leq u^\top x\}$, Tukey depth is therefore also called\emph{ halfspace depth} in the literature.
%\medskip

$T^*(\mathcal{X}^n)$ inherits many advantages of the univariate median. For example, it has an asymptotical breakdown point as high as $1/3$ under the centro-symmetric assumption. Unlike the aforementioned generalizations, $T^*(\mathcal{X}^n)$ is affine equivariant. That is,
\begin{eqnarray*}
    T^*(\mathbb{A}\mathcal{X}^n + \textbf{\emph{b}}) = \mathbb{A} T^*(\mathcal{X}^n) + \textbf{\emph{b}}
\end{eqnarray*}
for all full-rank $p\times p$ matrices $\mathbb{A}$ and $p$-vectors $\textbf{\emph{b}}$. Here $\mathbb{A}\mathcal{X}^n = \{\mathbb{A}X_1, \mathbb{A} X_2, \cdots, \mathbb{A}X_n\}$. Since data transformations, such as rotation and rescaling, are very common in practice, this property is very important and often expected for a multivariate location estimator.
\medskip

Unfortunately, $T^*(\mathcal{X}^n)$ is computationally intensive. Approximate algorithm has been developed by \cite{SR2000} though. This algorithm is lack of affine equivariance, nevertheless. While, on the other hand, efficient \emph{exact} algorithm exists, but is only feasible for \emph{bivariate} data \citep{RR1998, MRS2003}.
\medskip

Following DG92, define the $\alpha$th depth contour/region as
$D_\tau:=\{\x\in \mathcal{R}^p : D(\x,\mathcal{X}^n)\geq \tau\}$. To avoid unbounded as well as empty regions, hereafter $D_\tau$ is restricted to $0< \tau \leq \lambda^*(\mathcal{X}^n)$ \citep{LMM2014}.
Note that $T^*(\mathcal{X}^n)$ is in fact the average of all points in the most central (deepest) region, hereafter median region,
\begin{eqnarray*}
    \mathcal{M}(\mathcal{X}^n) = \left\{\x \in \mathcal{R}^p : D(\x, \mathcal{X}^n) = \lambda^*(\mathcal{X}^n) \right\}.
\end{eqnarray*}

Computing $T^*(\mathcal{X}^n)$ is in fact a special case of computing the depth regions. For a fixed depth value $\tau$, there have been some exact algorithms developed for computing the depth region \citep{PS2012a, PS2012b, RR1996}.
%
%Its computing issue can be transformed into that of computing the Tukey depth contours. For a depth contours with a fixed depth value $\tau$, there have been some efficient algorithms.
However, the maximum Tukey depth $\lambda^*(\mathcal{X}^n)$ is usually \emph{unknown} in advance. Hence, how to determine $\lambda^*(\mathcal{X}^n)$ is the most key issue for exactly computing $T^*(\mathcal{X}^n)$.
\medskip

To obtain this value, based on DG92's lower and upper bounds results, conventionally one needs to search over the interval $[\lceil n / (p + 1) \rceil / n, \lceil n / 2 \rceil / n]$ step by step as suggested and did in \cite{RR1998}, where $\lceil\cdot \rceil$ denotes the ceiling function. During this process, a lot of depth regions has to be computed. Usually, computing a single depth region is very time-consuming especially when $p$ is large. %It is usually very time-consuming.
Therefore, the computation would be greatly benefited in the practice by narrowing the search scope of the maximum Tukey depth, if possible.
\medskip

DG92 pointed out that the lower bound $\lceil n / (p + 1)\rceil / n$ of $D(x, \mathcal{X}^n)$ is attained for some special data sets in general position, and hence can not be further improved. Is there any room of improvement for the upper bound $\lceil n / 2 \rceil / n$ given in DG92?
\medskip

%In general, the exact maximum depth $\lambda^*(\mathcal{X}^n)$ is unknown. How to determine it is nontrivial. Hence, determining the %sharp bounds of the maximum depth would still

%It is almost impossible to to obtain it. The feasible way is to determine possible boundaries of $\lambda^*(\mathcal{X}^n)$. In the %literature, \cite{DG1992} have obtained a then DG gave an upper bound. But it is not sharp and could waste lots of step in %construction of depth contours. Here of course, we have to list some reference that did the extra steps. That is the references %which deal with HD contours.
%\medskip

The upper bound $\lceil n / 2 \rceil / n$ does not depend on the dimension $p$. Intuitively, this is not sensible for computing Tukey median when $p$ is increasing. Because for a data set with fixed $n$, an increasing proportion of observations will be pushed onto the boundary of the convex hull of the data points when $p$ increases \citep{HTF2008}. Hence, the maximum depth value may decrease. The upper bound given by DG fails to reflect this change in trend.
\medskip

In this paper, we provide a sharper upper bound ${\lfloor (n - p + 2) / 2 \rfloor}/{n}$ which may be much less than $\lceil n / 2 \rceil / n$ when $p$ is large relative to large $n$. For example, when $n = 5p$ with $p = 5$, the search range of the maximum Tukey depth can be reduced by \emph{more than one quarter}.
%\begin{eqnarray*}
%    (\lceil n / 2 \rceil - \lfloor (n - p + 2) / 2 \rfloor) / (\lceil n / 2 \rceil - \lceil n / (p + 1)\rceil) \ge 25\%.
%\end{eqnarray*}
%\cred{That is, the search scope of the maximum Tukey depth can be much narrowed.
%This is of great practical importance.
As for the choice of $n = 5p$,  a justification is given on p.326 of \cite{JP1995}; see also \cite{Zuo2004}.
\medskip

%
%Observe that for $p \ge 3$, the bound ${\lfloor (n - p + 2) / 2 \rfloor}/{n}$ is less than $1/2$ with respect to any given %$\mathcal{X}^n$ in general position. An interesting byproduct here is that `\emph{in general position}' and `\emph{halfspace %symmetry}' could not coexist for a data set in three or higher dimensions, even if it has an underlying central symmetric %distribution, e.g., normal distribution. These two notions are commonly assumed in the literature related to statistical depth %functions, and will be specified in the sequel.
\medskip

Since the computation of the depth of a single point is of less time complexity than that of computing a Tukey depth contour \citep{LZ2014a, DM2014}, and the depth of sample points is easily obtained as by-product of computing the depth contour, one may wonder if the deepest observation could serve as a Tukey median. This paper also considers this problem and provides a definite answer.
\medskip

We show that the sample observation can not lie in the interior of the median region. It can be included in the median region only if it is a vertex of the convex set.  Hence, unless median region is a singleton, a sample point can not serve as the Tukey median in general.
\medskip

To guarantee the uniqueness of Tukey's median for the population distributions, some symmetry assumption is commonly imposed( see, e.g., DG92, and \cite{ZS200a}(ZS00a)). The weakest version of symmetry is so-called `\emph{halfspace symmetry}'(see ZS00a). Another common assumption on underlying data when one deals with breakdown point robustness or data depth is `\emph{in general position}' (see DG92).\emph{ An interesting byproduct} of the paper is the conclusion: that `\emph{in general position}' and `\emph{halfspace symmetry}' could not coexist for a data set in three or higher dimensions, even if there is a central symmetric underlying distribution, e.g., normal distribution.
\medskip

The rest paper is organized as follows. Section 2 establishes a sharper upper bound for the halfspace depth of Tukey median.
Section 3 is devoted to provide a definite answer to the question: Can a sample point serve as Tukey median? Concluding remarks end the paper.

\vskip 0.1 in
\section{The upper bound for the Tukey depth}
\paragraph{}
\vskip 0.1 in \label{MainIdea}

In this section, we investigate the maximum Tukey depth. In the sequel we will always assume that data set is \emph{in general position}. A $p$-variate data set is called to be \emph{in general position} if there are no more than $p$ sample points in any $(p - 1)$-dimensional hyperplane. This assumption is often adopted in the literature when dealing with the data depth and breakdown point; see, e.g., \cite{DG1992, MLB2009}, among others.
\medskip

Observe that $n$ and the data cloud $\mathcal{X}^n$ are fixed in the following. For convenience, \emph{we drop the argument $\mathcal{X}^n$ from $T^*(\mathcal{X}^n)$, $\lambda^*(\mathcal{X}^n)$, $\mathcal{M}(\mathcal{X}^n)$ and $D(x, \mathcal{X}^n)$ if no confusion arises}. Besides, we introduce a few other notations as follows. Let $\kappa^* = n\lambda^*$. Obviously, $\kappa^*$ is a positive integer since the image of $P_n$ only can take a finite set of values $\{0, 1/n, 2/n, \cdots, 1\}$. For any $k$ points $\{\mathbf{x}_1, \mathbf{x}_2, \cdots, \mathbf{x}_k\}$, picking an arbitrary point, without loss of generality, say, $\mathbf{x}_k$,  we denote $\textbf{span}(\mathbf{x}_l|l = 1, 2, \cdots, k) = \{\sum_{l = 1}^{p-1} \lambda_l (\mathbf{x}_l - \mathbf{x}_k): \forall \lambda_i \in \mathcal{R}^1\}$ as
the subspace spanned by $\mathbf{x}_1 - \mathbf{x}_k, \mathbf{x}_2 - \mathbf{x}_k, \cdots, \mathbf{x}_{k-1} - \mathbf{x}_k$. (For $k = 1$, we assume $\textbf{span}(\mathbf{x}_l) = \{\mathbf{0}\}$.)
%\textcolor{blue}{Here again I want to raise this notation issue, this is not common definition for a spanned space, is there any advantages or reasons? Can one switch $x_i$ with $x_k$? Is your definition labeling(subscripts)-dependant? Or some other reference used this one?  Please provide the reference.}
\medskip

The following representation of $\mathcal{M}$ will be repeatedly utilized in the sequel. Following Theorem 4.2 in \cite{PS2011}, a finite number of direction vectors suffice for determining the sample Tukey regions, which include the median region as a special case. That is, let $\{\u_1^*, \u_2^*, \cdots, \u_m^*\}$ be all normal vectors of the hyperplanes, with each of which passing through $p$ observations and cutting off exactly $\kappa^* - 1$ observations, %{\color{green}(This statement is correct but is very difficult for readers to understand, e.g. what does ``cutting off'' means? Plus the following math expressions do not provide any further explanations or help, I would write like this: $\mathcal{M}(\mathcal{X}^n)=\bigcap_{j=1}^m H_j$, where $H_j$ is determined by $p$ observations and $n*P_n(H_j^C)=\kappa^*(\mathcal{X}^n) - 1$. Please modify this part. Original paper may much more clear than here since we take off some middle part from the original entire paper)}
then we have
\begin{eqnarray}
\label{RepreMXn}
    \mathcal{M} &=& \bigcap_{j=1}^m \left\{ \x: (\u_j^*)^\top \x \ge q_j \right\}, \quad \text{where}\\
    q_j &=& \inf \left\{t \in \mathcal{R}^1: P_n((\u_j^*)^\top X \leq t) \ge \lambda^* \right\}.\nonumber
\end{eqnarray}

The following theorem provides a sharp upper bound for the maximum Tukey depth.
\medskip

\textbf{Theorem 1}. Suppose $\mathcal{X}^n \subset \mathcal{R}^p$ ($n > p$) is in general position. We have
\begin{eqnarray*}
    \lambda^* \leq
    \left\{ \begin{array}{lcl}
    \frac{\lfloor (n - p + 2) / 2 \rfloor}{n}, &\quad & \textbf{dim}(\mathcal{M}) < p,\\
    & & \\
    \frac{\lfloor (n - p + 1) / 2 \rfloor}{n}, &\quad & \textbf{dim}(\mathcal{M}) = p,\\
    \end{array}
    \right.
\end{eqnarray*}
where $\textbf{dim}(\mathcal{M})$ denotes the affine dimension of $\mathcal{M}$.
\medskip

\begin{figure}[H]
\label{fig:IllForTh1}
\centering
	\subfigure[Illustration for Part (\textbf{I})]{
	\includegraphics[angle=0,width=2.9in]{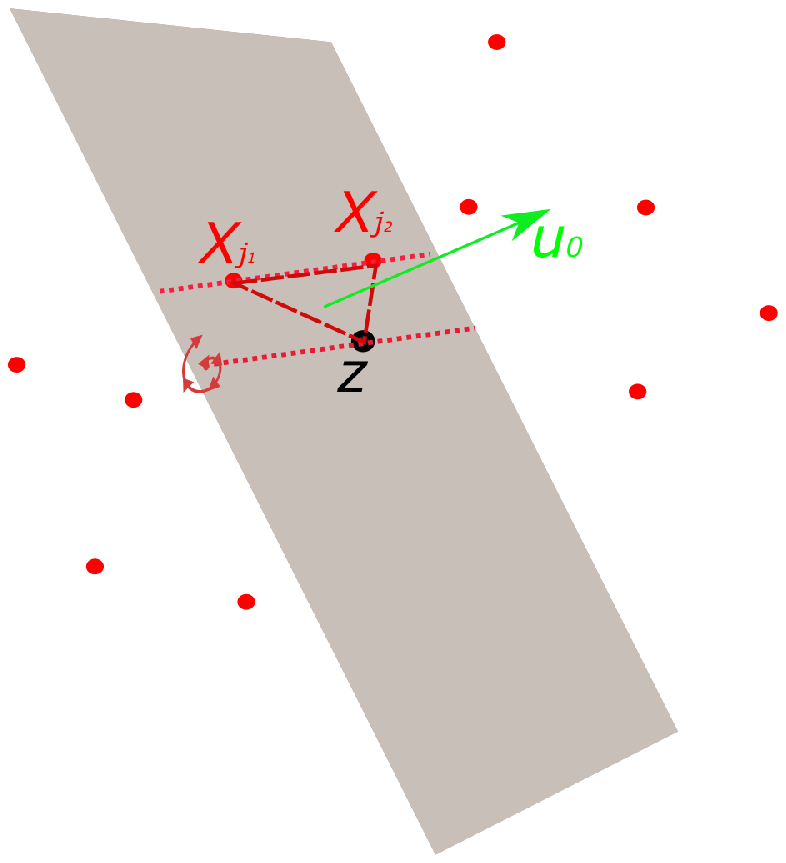}
	\label{fig:IllTh1001}
	}\quad
	\subfigure[Illustration for Part (\textbf{II})]{
	\includegraphics[angle=0,width=2.9in]{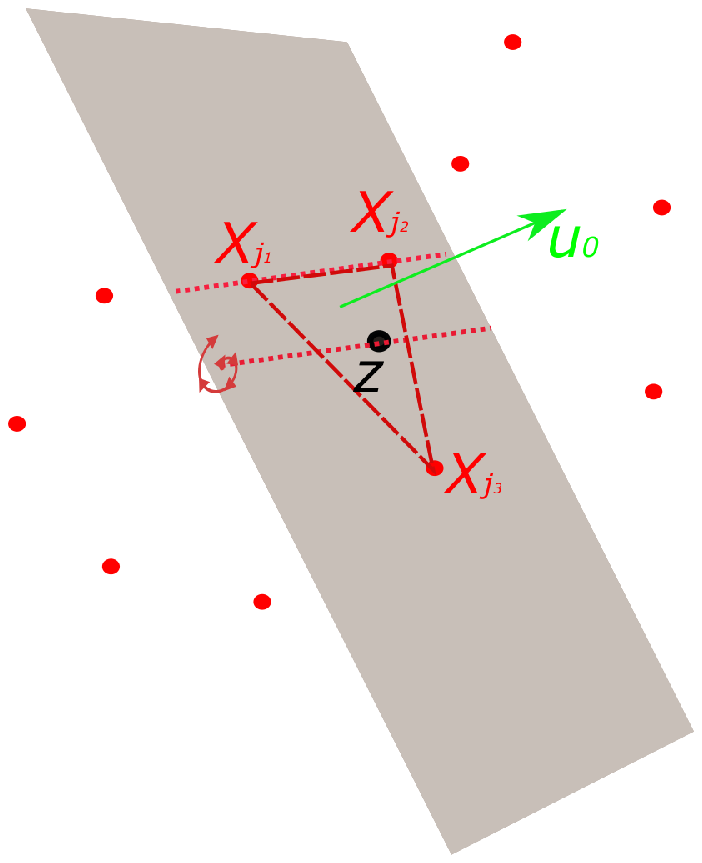}
	\label{fig:IllTh1002}
	}%\vspace{-0.8in}
\caption{Shown are illustrations for the main idea of the proofs of Part (\textbf{I}) (left) and (\textbf{II}) (right) of Theorem 1, where the small points denote the observations, the big point denote $\z$, and $\u_0$ the vector normal to the separating hyperplane, which passes through tow observations in the left case, but three observations in the right case.}
\end{figure}

Since the proof is trivial for $p=1$, we focus only on $p \ge 2$ in the sequel. For convenience, we present the long proof in two parts, i.e., (\textbf{I}) and (\textbf{II}). %Part (\textbf{I}) investigates the scenario $\textbf{dim}(\mathcal{M}(\mathcal{X}^n)) = p$, while Part (\textbf{II}) considers $\textbf{dim}(\mathcal{M}(\mathcal{X}^n)) < p$.
The basic idea is that: In Part (\textbf{I}), if $\mathcal{M}$ is of affine dimension $p$, i.e., has nonzero volume, one can always deviate (shift) the separating hyperplane around some points in $\mathcal{M}$ to get rid of $p-1$ sample points in it; See \ref{fig:IllTh1001}. In Part (\textbf{II}), as $\mathcal{M}$ is of affine dimension less than $p$, i.e., has volume zero, $T^*$ lies in a hyperplane through $p$ observations (not more because of the general position assumption), and thus one can not always cut off all of these $p$ observations because one point should remain in the lower mass halfspace; See \ref{fig:IllTh1002}.

\textbf{Proof of Theorem 1}. (\textbf{I}) $\textbf{dim}(\mathcal{M}) = p$: %Denote $\mathcal{M}^+ = \mathcal{M}(\mathcal{X}^n) \cap \{\x: \u^\top \x > \u^\top X_{i_1}\}$ and $\mathcal{M}^- = \mathcal{M}(\mathcal{X}^n) \cap \{\x: \u^\top \x < \u^\top X_{i_1}\}$. Then the convexity of $\mathcal{M}(\mathcal{X}^n)$ implies that at least one of $\mathcal{M}^+$ and $\mathcal{M}^-$ should be of affine dimension $p$. Without loss of generality, assume $\textbf{dim}(\mathcal{M}^+) = p$. Similarly, $\mathcal{M}^+$ can be further divided by the hyperplane passing through $p$ other observations into two parts, with at least one of them being of affine dimension $p$. Because the number of such combinations of $p$ observations is finite, i.e., ${n \choose p}$, there must exist one point $\z_1 \in \mathcal{M}(\mathcal{X}^n)$ such that it is \emph{not} contained by any hyperplane passing through $p$ observations. That is, $\{\z_1\} \cup \mathcal{X}^n$ is in general position.
Observe that $\mathcal{X}^n$ is a finite sample in general position, and $\mathcal{M}$ is a convex polytope of affine dimension $p$. Hence, it is easy to check that there exists a point $\z_1 \in \mathcal{M}$ such that $\{\z_1\} \cup \mathcal{X}^n$ is in general position.
\medskip

Using this, we claim that the hyperplane $\Pi_1$ passing through $\{\z_1, X_{j_1}, X_{j_2}, \cdots, X_{j_{p-1}}\}$ contains only $p-1$ observations of $\mathcal{X}^n$ for an any given set of $\{X_{j_1}, X_{j_2}, \cdots, X_{j_{p-1}}\}$. Let $\u_1$ be $\Pi_1$'s normal vector such that $P_n(u_1^\top X < \u_1^\top z_1) = \min\{P_n(\u_1^\top X < \u_1^\top \z_1), P_n(-\u_1^\top X < -\u_1^\top z_1)\}$. Obviously, $\u_1^\top \z_1 = \u_1^\top X_{j_1} = \cdots = \u_1^\top X_{j_{p-1}}$ and
%one can find a vector $\u_1$ satisfying that: (a) $\u_1^\top \z_1 = \u_1^\top X_{j_1} = \cdots = \u_1^\top X_{j_{p-1}}$ and (b) $P_n(u_1^\top X < \u_1^\top z_1) = \min\{P_n(\u_1^\top X < \u_1^\top \z_1), P_n(-\u_1^\top X < -\u_1^\top z_1)\}$.
%let $\vv_1 = \z_1 - X_{j_1}$. Under the condition of $\{\z_1\} \cup \mathcal{X}^n$ being in general position, there must exist a direction vector $\u_1$ determined uniquely by $\{\z_1, X_{j_1}, X_{j_2}, \cdots, X_{j_{p-1}}\}$ such that $\u_1^\top \z_1 = \u_1^\top X_{j_1} = \cdots = \u_1^\top X_{j_{p-1}}$ and $P_n(u_1^\top X < \u_1^\top z_1) = \min\{P_n(\u_1^\top X < \u_1^\top \z_1), P_n(-\u_1^\top X < -\u_1^\top z_1)\}$.
$P_n(\u_1^\top X < \u_1^\top \z_1) \leq \left\lfloor \frac{n - (p - 1)}{2} \right\rfloor \big/ n$.
\medskip

Let $\mathbf{V} = \textbf{span}(X_{j_k}| k = 1, 2, \cdots, p - 1)$. \emph{Now, we proceed to show that, if $\z_1 \notin \mathbf{V}$, it is possible to get rid of $X_{j_1}, X_{j_2}, \cdots, X_{j_{p-1}}$ from $\Pi_1$ through deviating it around $\z_1$}.
\medskip

Decompose $\vv_1 := \z_1 - X_{j_1} = \vv_1' + \vv_1''$, where $\vv_1' \in \mathbf{V}$ and $\vv_1''\in \mathbf{V}^\bot$. Here $\mathbf{V}^\bot$ denotes the orthogonal complement space of $\mathbf{V}$. Obviously, $\vv_1'' \neq 0$. (If not, $\vv_1' = \z_1 - X_{j_1} \in \mathbf{V}$, contradicting with that $\{\z_1\} \cup \mathcal{X}^n$ is in general position.) Similarly to \cite{DM2014}, let
\begin{eqnarray}
\label{VarEp}
    \varepsilon = \frac{1}{2} \min_{l \in \{1,\, 2,\, \cdots,\, n\} \setminus \{j_1,\, j_2,\, \cdots,\, j_{p-1}\}} \frac{|\u_1^\top (\z_1 - X_l)|}{|(\vv_1'')^\top (\z_1 - X_l)|},
\end{eqnarray}
and $\widetilde{\u}_1 = \u_1 - \varepsilon \vv_1''$. Here we define $|\u_1^\top (\z_1 - X_i)| / |(\vv_1'')^\top (\z_1 - X_i)| = +\infty$ if $(\vv_1'')^\top (\z_1 - X_i) = 0$. Then $ |\varepsilon(\vv_1'')^\top (\z_1 - X_i)| < |\u_1^\top (\z_1 - X_i)|$, and in turn
\begin{eqnarray}
\label{zVV}
    \text{sgn}(\widetilde{\u}_1^\top (\z_1 - X_i)) = \text{sgn}(\u_1^\top (\z_1 - X_i)), \quad i \in \{1,\, 2,\, \cdots,\, n\} \setminus \{j_1,\, j_2,\, \cdots,\, j_{p-1}\},
\end{eqnarray}
where sgn$(\cdot)$ denotes the sign function. On the other hand, (a) $\widetilde{\u}_1^\top (\z_1 - X_{j_1}) = -\varepsilon \|\vv_1''\|^2 < 0$, and (b) for each $k = 2, 3\cdots, p - 1$, $\widetilde{\u}_1^\top (\z_1 - X_{j_k}) = \widetilde{\u}_1^\top (\z_1 - X_{j_1} + (X_{j_1} - X_{j_k})) = -\varepsilon \|\vv_1''\|^2 < 0$. Then $\widetilde{\u}_1^\top \z_1 < \widetilde{\u}_1^\top X_{j_1}, \widetilde{\u}_1^\top X_{j_2}, \cdots, \widetilde{\u}_1^\top X_{j_{p-1}}$. These, together with \eqref{zVV}, lead to $P_n(\widetilde{\u}_1^\top X \leq \widetilde{\u}_1^\top \z_1) = P_n(\u_1^\top X < \u_1^\top \z_1)$. Hence,
\begin{eqnarray*}
    \lambda^* = D(\z_1) \leq P_n(\widetilde{\u}_1^\top X \leq \widetilde{\u}_1^\top \z_1) \leq \left\lfloor\frac{n - (p - 1)}{2} \right\rfloor \Big/ n.
\end{eqnarray*}

(\textbf{II}) $\textbf{dim}(\mathcal{M}) < p$: Relying on \eqref{RepreMXn}, we claim that there $\exists k \in \{1, 2, \cdots, m\}$ such that $(\u_k^*)^\top T^*(\mathcal{X}^n) = q_k$. Using this, a simple derivation leads to $\mathcal{M} \subset \Pi_2 := \{z \in \mathcal{R}^p : (\u_k^*)^\top \z = q_k\}$. Obviously, $\Pi_2$ should pass through $p$ observations, say $\{X_{k_1}, X_{k_2}, \cdots, X_{k_p}\}$, by the definition of $\u_2 := \u_k^*$ and $q_k$, %For a given $z_2 \in \mathcal{M}(\mathcal{X}^n)$, let $u_2$ be the normal vector of $\Pi_0$ satisfying that $P_n(u_2^\top X < u_2^\top z_2) = \min\{P_n(u_2^\top X < u_2^\top z_2), P_n(-u_2^\top X < -u_2^\top z_2)\}$.
and $P_n(\u_2^\top X < \u_2^\top \z_2) \leq \lfloor \frac{n - p}{2} \rfloor \big/ n$ for any $\z_2 \in \mathcal{M}$.
\medskip

Without confusion, assume $X_{k_p} = 0$ by the affine equivariance of the Tukey depth. Let $\mathbf{W}_i = \textbf{span}(X_{k_l}| l \in \{1, 2, \cdots, p\} \setminus \{i\})$ for $i = 1, 2, \cdots, p$. We now show that $\z_2$ does \emph{not lie simultaneously in} all affine subspaces $\mathbf{W}_i$'s defined by facets of a $(p-1)$-dimensional simplex formed by $\{X_{k_1}, X_{k_2}, \cdots, X_{k_p}\}$.
\medskip

\emph{For $p = 2$}: The proof is trivial because $\z_2 - X_{k_1}$, $\z_2 - X_{k_2}$ could not be $0$ simultaneously.
\medskip

\emph{For $p > 2$}: By observing that if $\z_2 = 0 ( = X_{k_p})$,  $\z_2 \notin \mathbf{W}_p$ because $\mathcal{X}^n$ is in general position, we only focus on the case $\z_2 \neq 0$ in what follows. %We will prove that the assumption that $\z_2 \in \mathbf{W}_i$ for all $i = 1, 2, \cdots, p$ leads to a contradiction.
\medskip

Suppose $\z_2 \in \mathbf{W}_i$ for all $i = 1, 2, \cdots, p$, then there must exist some constants such that:
\begin{eqnarray}
\label{subspace}
    \left\{ \begin{array}{cclcl}
    \z_2 &=& \quad a_{12}X_{k_2} + a_{13}X_{k_3} + \cdots + a_{1, p - 1}X_{k_{p-1}} &\quad & (z_2 \in \mathbf{W}_1)\\
    \z_2 &=& a_{21}X_{k_1} + \quad a_{23}X_{k_3} + \cdots + a_{2, p - 1}X_{k_{p-1}} &\quad & (z_2 \in \mathbf{W}_2)\\
    & \vdots & &\quad & \\
    \z_2 &=& a_{p-1,1}X_{k_1} + a_{p-1, 2}X_{k_2} + \cdots + a_{p-1, p - 2}X_{k_{p-2}} \quad &\quad & (z_2 \in \mathbf{W}_{p-1})\\
    \z_2 &=& b_{1}X_{k_1} + b_{2}X_{k_2} + \cdots + b_{p-1}X_{k_{p-1}} &\quad & (z_2 \in \mathbf{W}_p).
    \end{array}
    \right.
\end{eqnarray}
Let $(w_1, w_2, \cdots, w_{p-1})^\top$ be a solution to
\begin{eqnarray*}
    \left( \begin{array}{ccccc}
    0& a_{21}& a_{31} & \cdots & a_{p-1,1}\\
    a_{12}& 0& a_{32} & \cdots & a_{p-1,2}\\
    a_{13}& a_{23}& 0& \cdots & a_{p-1,3}\\
    \vdots& \vdots& \vdots & \ddots & \vdots\\
    a_{1,p-1}& a_{2,p-1}& a_{3,p-1} & \cdots & 0
    \end{array}\right)
    \left( \begin{array}{c}
    w_1\\
    w_2\\
    w_3\\
    \vdots\\
    w_{p-1}
    \end{array}\right)
    =
    \left( \begin{array}{c}
    b_1\\
    b_2\\
    b_3\\
    \vdots\\
    b_{p-1}
    \end{array}\right).
\end{eqnarray*}
Then it is easy to check that $\sum_{i=1}^{p-1} w_i \z_2 - \z_2 = 0$. Hence, $\sum_{i=1}^{p-1} w_i = 1$ due to $\z_2 \neq 0$. On the other hand, $\z_2 \neq 0$ implies $(b_1, b_2, \cdots, b_{p-1}) \neq 0$. Without confusion, assume $b_{i_0} \neq 0$ for $1 \leq i_0 \leq p-1$. Then by letting `the $p$-th equation' $-$ `$\sum_{i=1}^{p-1}w_i$ $\times$ the $i_0$-th equation' of \eqref{subspace}, we can obtain
\begin{eqnarray*}
    X_{k_p} = 0 = \z_2 - \sum_{i=1}^{p-1} w_i\z_2 = \sum_{l=1}^{p-1} b_l X_{k_l} - \sum_{i=1}^{p-1} w_i\sum_{l\neq i_0}^{p-1} a_{i_0,l} X_{k_l} =: \sum_{l \neq i_0}^{p - 1} \widetilde{b}_l X_{k_l} + b_{i_0} X_{k_{i_0}} \in \mathbf{W}_p.
\end{eqnarray*}
%where $\widetilde{b}_l = b_l - \sum_{i=1}^{p-1} w_i a_{i_0, l}$ for $1 \leq l \leq p - 1$ but $l \neq i_0$.
This clearly \emph{contradicts with} the fact $X_{k_p} \notin \mathbf{W}_p$ because $\mathcal{X}^n$ is in general position.
\medskip

Suppose $\z_2 \notin \mathbf{W}_p$ without confusion. Then similar to Part (\textbf{I}), it is possible to get rid of $X_{k_1}, X_{k_2}, \cdots, X_{k_{p-1}}$ from the separating hyperplane $\Pi_2$ through rotating it around $\z_2$, and hence $\lambda^* = D(\z_2) \leq P_n(\u_2^\top X < \u_2^\top \z_2) + 1 / n \leq \left\lfloor\frac{n - p + 2}{2} \right\rfloor \big/ n$.
\medskip

This completes the proof of Theorem 1. \hfill $\Box$

\begin{figure}[H]
\centering
\includegraphics[angle=0,width=4.0in]{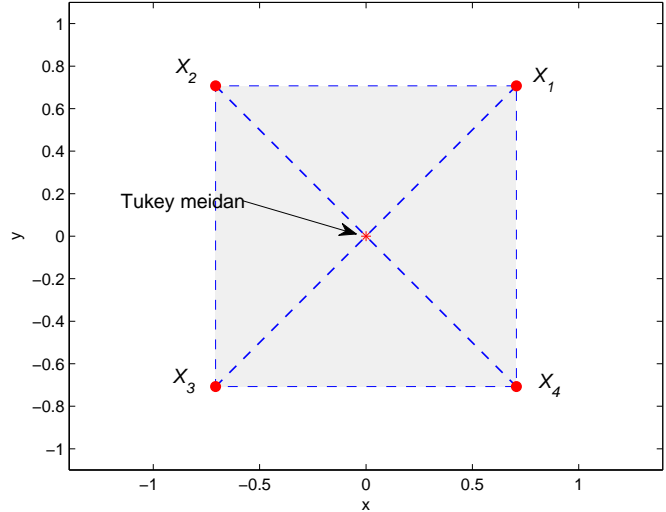}
\caption{The maximum Tukey depth in this example is $\frac{\lfloor (4 - 2 + 2) / 2 \rfloor}{4} = 1/2$.}
\label{fig:HM2D003}
\end{figure}

\noindent
\textbf{Remarks}\\[1ex]
\textbf{(i)} The upper bound given in Theorem 1 is \emph{attained} if the data set is strategically chosen; see Figure \ref{fig:HM2D003} for $p = 2$ and Figure \ref{fig:TM3D1} for $p = 3$, for example. Hence, this bound is \emph{sharp}, and \emph{can not be further improved if the data are in general position}. The upper bound $\lceil n / 2 \rceil/n$ for the maximum Tukey depth given in DG92 clearly is\emph{ not} sharp for $p \ge 3$ (see Figure \ref{fig:TM3D1}).

%(i) The upper bound given in Theorem 1 is \emph{attained} if the data set is strategically chosen; see Figure \ref{fig:HM2D003}, for %example. Hence, this bound is \emph{sharp}, and \emph{can not be further improved if the data are in general position}. The upper %bound $\lceil n / 2 \rceil/n$ for the maximum Tukey depth given in DG92 clearly is\emph{ not} sharp.

\noindent
\textbf{(ii)}
Since the upper bound is smaller than $\lceil n / 2 \rceil/n$ for $p\ge 3$, an interesting question is raised here about the multivariate symmetry, which is closely related to the multidimensional median.
\vskip 3mm

In the literature, various notions of multivariate symmetry already exist, e.g., central symmetry, angular symmetry and halfspace symmetry. A random vector $X\in \mathcal{R}^p$ is said to have a distribution centrally symmetric about $\theta$ if $X - \theta$ and $\theta - X$ are equal in distribution. The angular symmetry, introduced in \cite{Liu1990}, is broader than the central symmetry. A random vector $X \in \mathcal{R}^p$ has a distribution angularly symmetric about $\theta$ if $(X - \theta) / \|X - \theta\|$ has a centrally symmetric distribution about the origin. Among three multivariate symmetry notions, the halfspace symmetry is the most broadening. As introduced by \cite{ZS200b}, a random vector $X$ has a distribution halfspace symmetric about $\theta$ if
\begin{eqnarray}
\label{Hsymm}
    P(X \in \mathcal{H}_\theta) \ge 1/2, \text{ each closed halfspace } \mathcal{H}_\theta \text{ with } \theta \text{ on the boundary.}
\end{eqnarray}

It is readily to see that any one-dimensional data set is always halfspace symmetric about its median. Whether this property still holds in spaces with dimension $p>1$ is not clear without Theorem 1.
 \medskip

As a byproduct, Theorem 1 actually provides a negative answer to this question when $\mathcal{X}^n$ is in general position. By Theorem 1, for $p \ge 3$,
the maximum Tukey depth is less than $1/2$ with respect to any given $\mathcal{X}^n$ in general position,
 %It avoids the definition in \eqref{Hsymm}. Hence,
%the empirical distribution of
therefore $\mathcal{X}^n$ could not be halfspace symmetric (and consequently not central or angular symmetric) about a point $\theta$ when $\mathcal{X}^n$ is in general position, even if $\mathcal{X}^n$ are generated from a central symmetric distribution, e.g., normal distribution. That is, the following proposition holds.
\medskip

\textbf{Proposition 1}. \emph{In general position }and \emph{halfspace symmetry} could not coexist for $\mathcal{X}^n$ in $R^p$ with $p \ge 3$.
\medskip
%
%we have} that:
%\begin{enumerate}
%  \item[(\textbf{a})] If $\mathcal{X}^n$ is in general position, then  $\mathcal{X}^n$ could \emph{not} be halfspace symmetric (and %consequently not central or angular symmetric) about a point $\theta$.
%
%  \item[(\textbf{b})] If $\mathcal{X}^n$ are central/angular/halfspace symmetric about a point $\theta$, then $\mathcal{X}^n$ could %\emph{not} be in general position.
%\end{enumerate}

\noindent
\textbf{(iii)}
For a given data set in general position, although one does not know what is the median region nor its dimension, Theorem 1 still provides a useful guide: i.e. the maximum depth is less than $\lfloor(n-p+2)/2\rfloor/n$.

%\cred{(iii) A question raised here is how to determine $\textbf{dim}(\mathcal{M}(\mathcal{X}^n))$? How many cases there could be? Fortunately, there are only two cases: equal to $p$ or $0$. That is, $\mathcal{M}(\mathcal{X}^n)$ is either a complete subspace with dimension $p$ or a singlton (please see Theorem 3).}
%\medskip
%\textcolor{magenta}{Claim: When $\textbf{dim}(\mathcal{M})<p$, $\mathcal{M}$ actually is a singleton. Could you prove it or disprove it?}

\vskip 0.1 in
\section{Can the deepest sample point serve as the Tukey median?}
\paragraph{}
\vskip 0.1 in \label{MainIdea}

It is well known in the literature that computing the depth of a single point is of less time complexity than computing a Tukey depth region. Hence, it would save a lot of effort in computing and searching Tukey median if a single deepest point could serve the role. From last section, we see that $\mathcal{M}$ can be a singleton. Could a deepest sample point serve as Tukey median?  Results following will answer the question.
\medskip

\textbf{Theorem 2}. Suppose $\mathcal{X}^n$ is in general position. If there exists an observation $X_l$ ($1 \leq l \leq n$) such that $D(X_l) = \lambda^*$, then $X_l$ must be a vertex of $\mathcal{M}$.
\medskip

\textbf{Proof}. Denote $d := \textbf{dim}(\mathcal{M})$. Clearly, $0 \leq d \leq p$. When $d = 0$, the proof is trivial because $\mathcal{M}$ is a singleton. In the sequel we only focus on the case $d = p$. (For $0 < d < p$, the proof is similar.) Since the proof is too long, we divide it into two parts.
\medskip

\textbf{(S1)}. \emph{We first show that $X_l$ could not lie in the inner of $\mathcal{M}$}. By the compactness of $\{\z \in \mathcal{R}^p: \|\z\| = 1\}$ and the fact that the image of $P_n$ takes a finite set of values, one can show that there must exist a $\u_0$ such that
\begin{eqnarray*}
  P_n(\u_0^\top X \leq \u_0 X_l) = \lambda^*,
\end{eqnarray*}
i.e., there exist $i_1,\, i_2,\, \cdots,\, i_{n - \kappa^*}$ satisfying that $\u_0 X_l < \u_0 X_{i_1} \leq \u_0 X_{i_2} \leq \cdots \leq \u_0 X_{i_{n - \kappa^*}}$.
\medskip

Let $\mathcal{H}(X_l, \u_0) = \{\x \in \mathcal{R}^p | \u_0^\top \x < \u_0^\top X_l\}$. Clearly, $\mathcal{H}(X_l, \u_0) \cap \mathcal{M} = \emptyset$. If not, there $\exists \z_0 \in \mathcal{M} \cap \mathcal{H}(X_l, \u_0)$ such that $\u_0^\top \z_0 < \u_0^\top X_l$, which implies $D(\z_0) \leq P_n(\u_0^\top X \leq \u_0^\top \z_0) < \lambda^*$. This contradicts with $\z_0 \in \mathcal{M}$. Hence, $X_l$ should be on a $(p - 1)$-dimensional facet $\mathcal{F}$ of $\mathcal{M}$.
\medskip

%With loss of generality, assume that $u_0 \perp \mathcal{F}$.
\textbf{(S2)}. \emph{In this part, we further show that $X_l$ should be on a $(p - 2)$-dimensional facet of $\mathcal{F}$}. Suppose $\mathcal{F}$ lies on the hyperplane $\Pi_3$ passing through $\{X_l\} \cup \{X_{j_1},\, X_{j_2},\, \cdots,\, X_{j_{p - 1}}\}$ and contains $m$ vertices $\{\mathbf{\vv}_s\}_{s = 1}^m =: \mathcal{V}$ of $\mathcal{M}$. Obviously, $\mathcal{F}$ is convex, because it is an intersection of halfspaces. For convenience, let $\mathbf{S} = \textbf{span}(X_{j_k}| k = 1, \cdots, p - 1)$ and $\mathbf{S}^\bot$ be its orthogonal complement space.
\medskip

\begin{figure}[H]
\centering
\includegraphics[angle=0,width=3.0in]{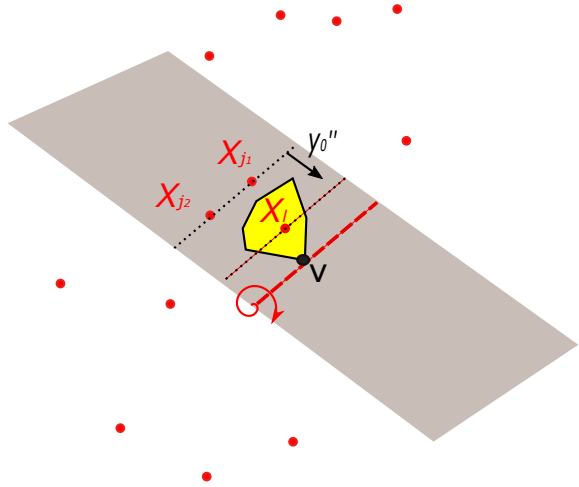}
\caption{Shown is a 3-dimensional illustration for Theorem 2 if $X_l$ is an inner point of a 2-dimensional facet $\mathcal{F}$, i.e., the polygon on the separating hyperplane, of $\mathcal{M}$.}
\label{fig:Th2}
\end{figure}

\emph{If $X_l$ lies in the inner of $\mathcal{F}$}, then (i) $\u_0$ should be normal to $\mathcal{F}$, and (ii) there $\exists \textbf{v} \in \mathcal{V}$ satisfying
\begin{eqnarray}
\label{OBV001}
    (\y_0'')^\top (\mathbf{v} - X_{j_k}) > \|\y_0''\|^2, \quad \text{for all } k = 1,\, \cdots,\, p - 1,
\end{eqnarray}
where $\y_0'' = \y_0 - \y_0' \in \mathbf{S}^\bot$ with $\y_0 = X_l - X_{j_1}$ and $\y_0' \in \mathbf{S}$. (\emph{If not}, for each $\mathbf{v}_s \in \mathcal{V}$, there $\exists \mathbf{x} \in \{X_{j_1},\, X_{j_2},\, \cdots,\, X_{j_{p - 1}}\}$ such that $\|\y_0''\|^2 \ge (\y_0'')^\top (\mathbf{v} - \mathbf{x}) = (\y_0'')^\top (\mathbf{v} - X_{j_1} + X_{j_1} - \mathbf{x}) = (\y_0'')^\top (\mathbf{v} - X_{j_1})$, which further implies, for $\forall \x = \sum_{s = 1}^m \lambda_s \mathbf{v}_s \in \mathcal{F}$,
\begin{eqnarray*}
    (\y_0'')^\top \left(\sum_{s = 1}^m \lambda_s \mathbf{v}_s - X_{j_1}\right) = \sum_{s = 1}^m \lambda_s \left((\y_0'')^\top (\mathbf{v}_s - X_{j_1}) \right) \leq \|\y_0''\|^2,
\end{eqnarray*}
where $\lambda_s \ge 0$ and $\sum_{s = 1}^m \lambda_s = 1$. \emph{This is impossible because $(\y_0'')^\top (X_l - X_{j_1}) = \|\y_0''\|^2$ and $X_l$ is an inner point of $\mathcal{F}$}.) That is, $\textbf{v}$ is farther away from $\mathbf{S}$ than $X_l$ along the same direction $\y_0''$. Using this and the fact $\|\y_0''\|^2 > 0$ due to the general position assumption of $\mathcal{X}^n$, a similar proof to Part (\textbf{I}) of Theorem 1 can show that it is possible to get rid of all $p$ points $X_l, X_{j_1},\, X_{j_2},\, \cdots,\, X_{j_{p - 1}}$ from $\Pi_3$ through deviating it around $\textbf{v}$; See Figure \ref{fig:Th2} for a 3-dimensional illustration. As a result, $D(\mathbf{v}) < P_n(\u_0^\top X \leq \u_0^\top X_l) = D(X_l)$, \emph{contradicting with $\mathbf{v} \in \mathcal{M}$}. Hence, $X_l$ should be on a $(p - 2)$-dimensional facet of $\mathcal{M}$.
\medskip

Similar to (\textbf{S2}), one can always obtain a contradiction if $X_l$ lies in the inner of a $k$-dimensional facet of $\mathcal{M}$ for $0 < k \leq p - 2$. This completes the proof. \hfill $\Box$
\medskip

\begin{figure}[H]
\centering
\includegraphics[angle=0,width=4.0in]{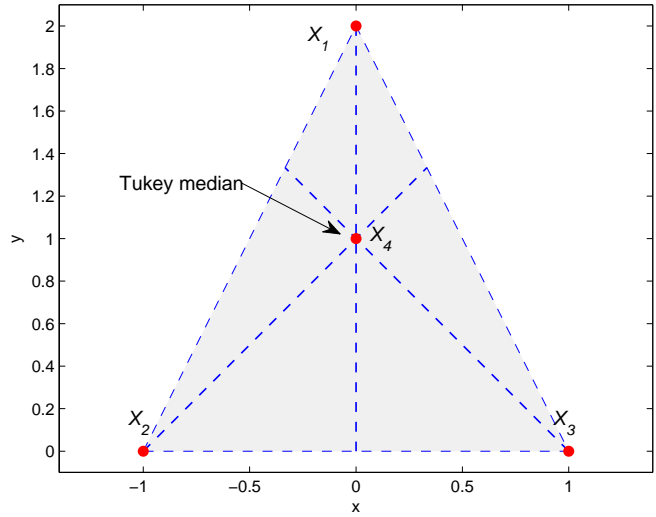}
\caption{Shown is an example such that $\mathcal{M}$ contains only a single sample point, namely, $X_4$. Hence, $X_4$ can serve as the Tukey median.}
\label{fig:Single}
\end{figure}

Theorem 2 indicates that the sample point could not lie in the interior of the median region $\mathcal{M}$. Hence, unless $\mathcal{M}$ contains only a single sample point (see Figure \ref{fig:Single} for an illustration), the sample point can not be used as the Tukey median which is defined to be the average of all point in $\mathcal{M}$.

\begin{figure}[H]
\centering
\includegraphics[angle=0,width=4.0in]{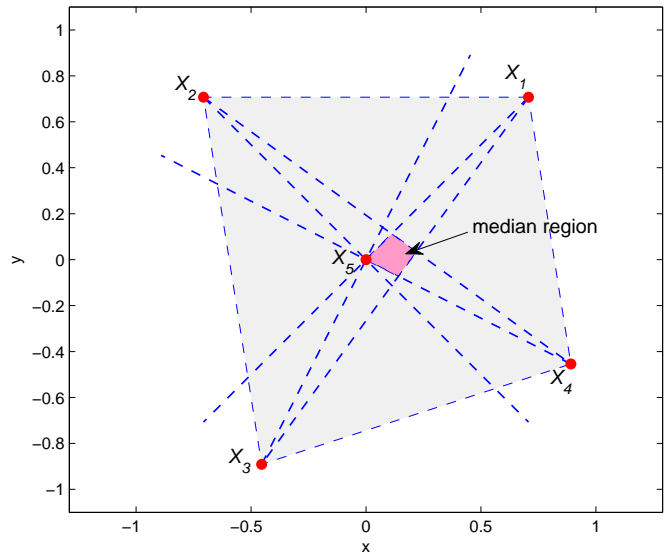}
\caption{Shown is an example such that $\mathcal{M}$ is of affine dimension $p=2$.}
\label{fig:GreatThan0}
\end{figure}

\begin{figure}[H]
\begin{center}
	\subfigure[Shown is an example such that $\mathcal{M}$ is singleton.]{
	\includegraphics[angle=0,width=2.8in]{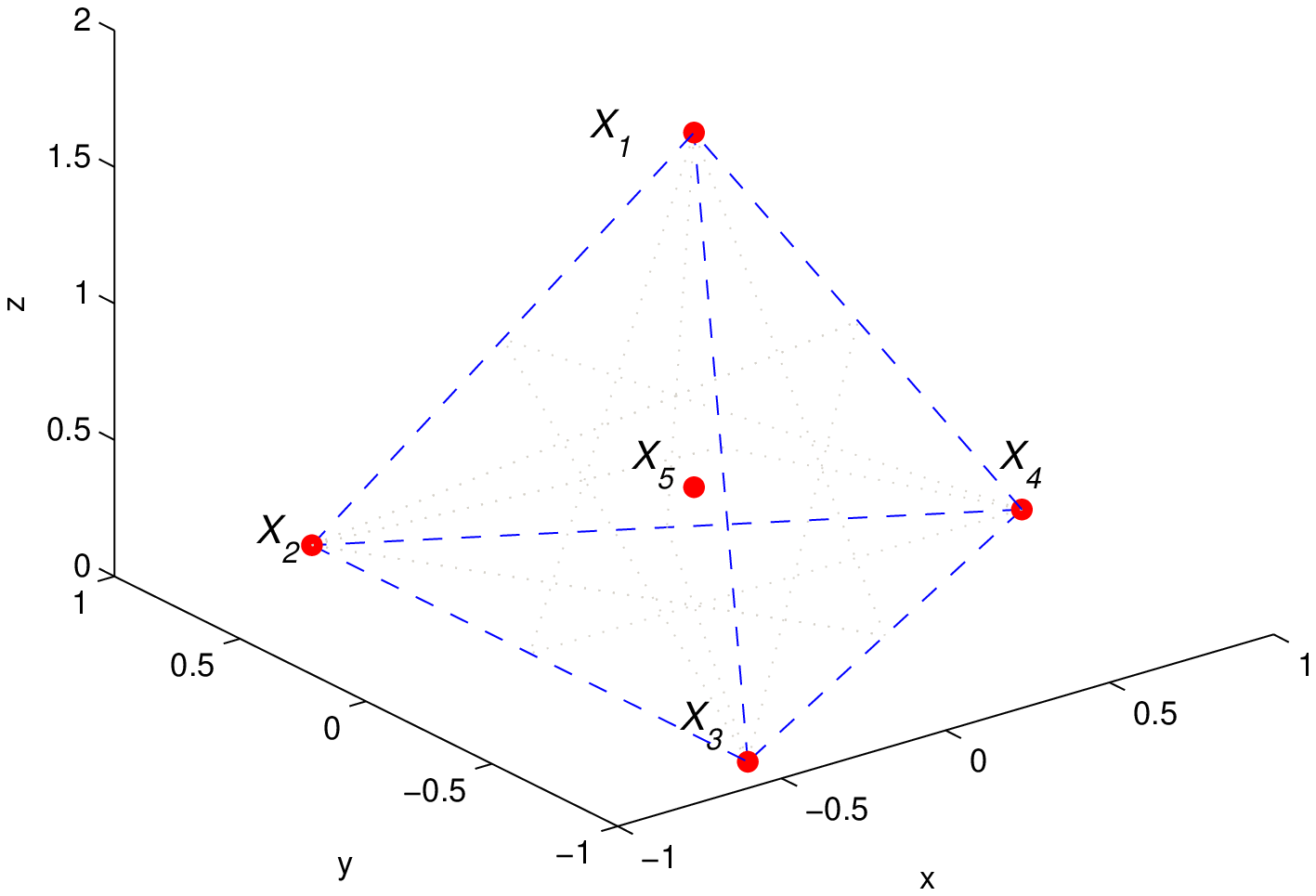}
	\label{fig:TM3D1}
	}\quad
	\subfigure[Shown is an example such that $\mathcal{M}$ is of affine dimension $p=3$.]{
	\includegraphics[angle=0,width=2.8in]{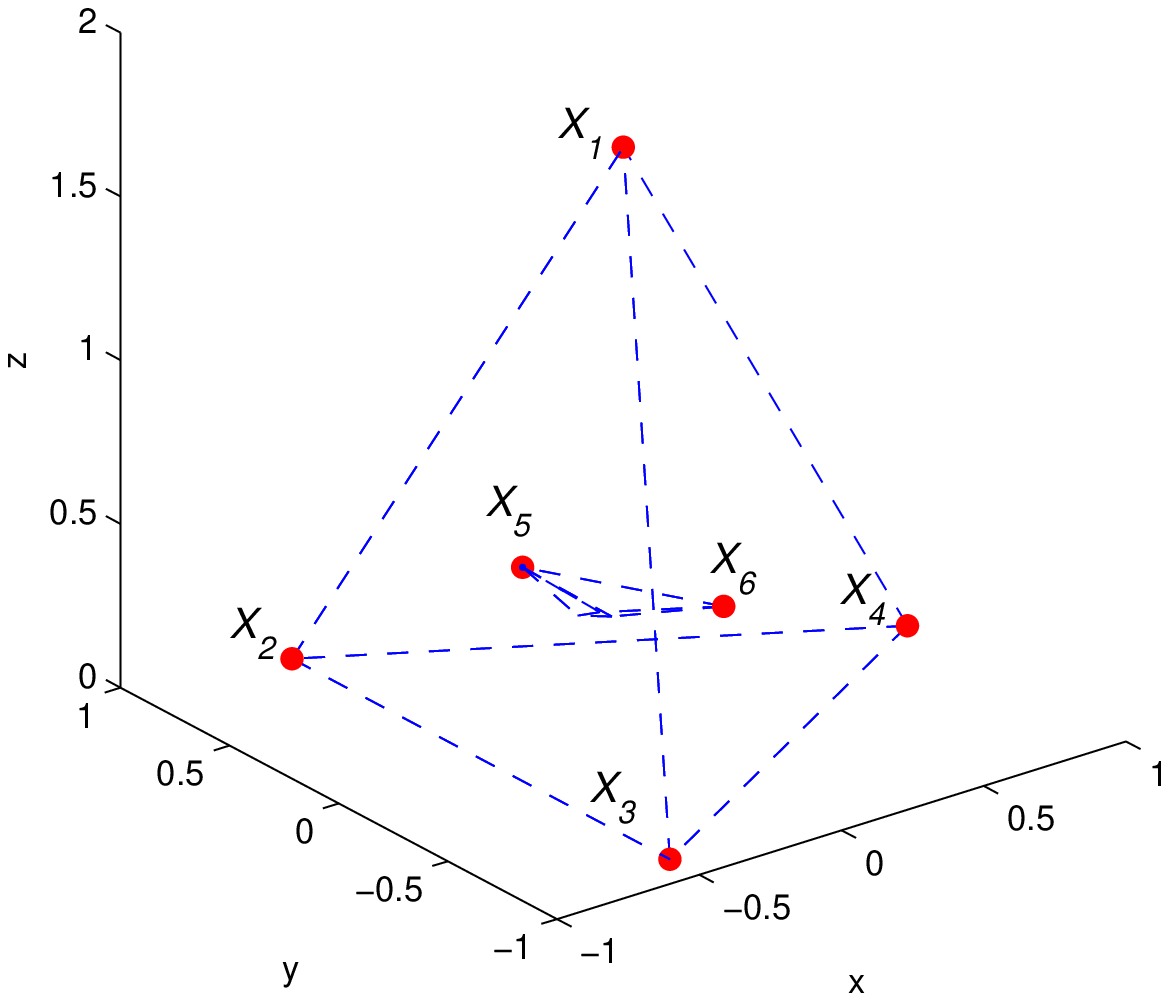}
	\label{fig:TM3D2}
	}
\caption{Shown are examples of $\mathcal{M}$ when $p = 3$.}
\label{fig:TM3DExam}
\end{center}
\end{figure}

Unfortunately, the latter scenario is very much possible in practice in the sense that $\mathcal{M}$ contains more than a single point. As one can see from Figure \ref{fig:GreatThan0}, although $X_5$ is one of the points maximizing the Tukey depth, it can not be used as the Tukey median for the sake of affine equivariance and because the average of $\mathcal{M}$ should be in the interior of $\mathcal{M}$, while $X_5$ is on the boundary. (Two 3-dimensional examples of both scenarios are shown in Figure \ref{fig:TM3DExam}.)
\medskip

So when $\mathcal{M}$ is a singleton? The following theorem partially answers this question.
\medskip

\textbf{Theorem 3}. Suppose $\mathcal{X}^n$ is in general position. Then when $\lambda^* = \frac{\lfloor (n - p + 2) / 2 \rfloor}{n}$ with $n > p \ge 2$, $\mathcal{M}$ contains only a single point.
\medskip

\textbf{Proof}. In the sequel we will show that if $\textbf{dim}(\mathcal{M}) > 0$, it would lead to a contradiction under the current assumptions. We focus only on the scenario $\textbf{dim}(\mathcal{M}) = 1$. The proof of the rest cases follows a similar fashion.
\medskip

When $\textbf{dim}(\mathcal{M}) = 1$, $\mathcal{M}$ is in fact a line segment. Denote $\x_1$, $\x_2$ to be its two endpoints. The compactness of $\mathcal{M}$ implies that, among $\{\u_1^*, \u_2^*, \cdots, \u_m^*\}$, there must exist a $\u_j^*$ such that
\begin{eqnarray}
\label{ell}
    (\u_j^*)^\top \x_1 = q_j \text{ and }(\u_j^*)^\top \x_2 > q_j.
\end{eqnarray}
If not, all points $\x$ in the line $\ell$ that passes through $\x_1$ and $\x_2$ satisfy that $(\u_k^*)^\top \x = q_k$, $k = 1, 2, \cdots, m$. This implies $\ell \subset \mathcal{M}$, contradicting with the boundedness of $\mathcal{M}$. On the other hand, by \eqref{RepreMXn}, there must exist $p$ observations, say $X_{i_1}, X_{i_2}, \cdots, X_{i_p}$, satisfying that $(\u^*)^\top X_{i_1} = (\u^*)^\top X_{i_2} = \cdots = (\u^*)^\top X_{i_p} = (\u^*)^\top \x_1$. (If $\x_1$ is a sample point, assume $X_{i_1} = \x_1$.) This, combined with \eqref{ell}, easily leads to
\begin{eqnarray*}
    n P_n((-\u_j^*)^\top X \leq (-\u_j^*)^\top \x_2) \leq \left\lfloor \frac{n - p}{2} \right\rfloor.
\end{eqnarray*}
This contradicts with the fact $\x_2 \in \mathcal{M}$. \hfill $\Box$
\medskip

%\textbf{Remark}
%\textcolor{magenta}{Claim in the end of last section (if it is true) is slightly broader than this result so we should replace this with the claimed result which I believe is true}.

\vskip 0.1 in
\section{Concluding remarks }
\paragraph{}
\vskip 0.1 in \label{Conclusion}

%\textcolor{magenta}{Here we need to at least report the time (whatever long it is) we used (with the fastest algorithm) to search the HM in $p=5$ and $n=5p$, then predict the time needed for the traditional approach and at the same time point out that single sample point could not serve as the halfspace median, this will make this manuscript more publishable. we can post the data and the median on our website if readers are interested}.
%\medskip

%\begin{figure}[H]
%\begin{center}
%\begin{minipage}{130mm}
%	\subfigure[Case 1.]{
%	\includegraphics[angle=0,width=2.5in]{TM3D1.eps}
%	\label{fig:ddplotXvsNorm}
%	}\quad
%	\subfigure[Case 2.]{
%	\includegraphics[angle=0,width=2.5in]{TM3D2.eps}
%	\label{fig:ddplotNormXSigm1}
%	}
%\caption{Shown is an example in dimension $p=3$.}
%\label{fig:DDplot}
%\end{minipage}
%\end{center}
%\end{figure}

%\begin{figure}[H]
%\centering
%\includegraphics[angle=0,width=4.0in]{TM3D1.eps}
%\caption{Shown is an example such that $\mathcal{M}(\mathcal{X}^n)$ is singleton when $p = 3$.}
%\label{fig:3dcase1}
%\end{figure}

%\begin{figure}[H]
%\centering
%\includegraphics[angle=0,width=4.0in]{TM3D2.eps}
%\caption{Shown is an example such that $\mathcal{M}(\mathcal{X}^n)$ is of affine dimension $p=3$.}
%\label{fig:3dcase2}
%\end{figure}

In the computing of Tukey's halfspace median, the lower and upper bounds given in DG92 on the maximum
halfspace depth are employed in the literature. The lower bound is sharp, but the upper bound is not in general (as shown in this paper), which could cause lots of unnecessary efforts in the searching of depth contours/regions. Computing of a singe depth contour can cost lots of time and effort. By providing a sharper and sharpest upper bound could save a lot of resource in practices, which is exactly achieved in this manuscript. Furthermore, we provide answers to the questions ``Can a single deepest sample point serve as Tukey's halfspace median? If yes, in what kind of situation?".
%In this paper, we studied the sample Tukey median. A sharper upper bound than the existing one is provided for the maximum Tukey depth. Data examples show that this bound could not be further improved under the in general position assumption. Two theorems are also provided for investigating the relationship between the deepest sample point and the Tukey median. %Note the maximum Tukey depth is larger than the depth of the deepest sample point.
\medskip

%In the literature, the Tukey median is known as its high breakdown point robustness. It turns out that its breakdown point may as high as $1/3$. Nevertheless, this value is only an asymptotic result. For a fixed $n$ (that is finite sample case), how large is this value is still open.
%Since the maximum Tukey depth is closely related to the finite sample breakdown point of the Tukey median, we hope that results of this paper can be helpful to answer the question.
 %is worthy of further study. This is our ongoing aim in next.

Results established here are not only interesting themselves theoretically but useful practically as well. %They are novel for Tukey's halfspace median, although similar derivations have been conducted by \cite{ZCH2004} and \cite{Mul2013} for breakdown point of some other medians.
Furthermore, observe that the finite sample breakdown point (FSBP) of Tukey's halfspace median $T^*$ is closely related to the maximum Tukey depth \citep{DG1992}. We anticipate that they will be extremely helpful for establishing the FSBP of $T^*$ and revealing its impact from the dimensionality $p$. This is still an open problem up to this point, although similar discussions have been conducted for some other multivariate location estimators \citep{ZCH2004, Mul2013}.

\vskip 0.1 in
\section*{Acknowledgements}
\paragraph{}
\vskip 0.1 in \label{Conclusion}

The research of the first two authors is supported by National Natural Science Foundation of China (Grant No.11461029, 61263014, 61563018), NSF of Jiangxi Province (No.20142BAB211014, 20143ACB21012,  20132BAB201011, 20151BAB211016), and the Key Science Fund Project of Jiangxi provincial education department (No.GJJ150439, KJLD13033, KJLD14034).
\medskip

We thanks the Chief-in-Editor Professor M\"uller, C., the AE, and two anonymous reviewers for their careful reading and insightful comments, which led to many improvements in this paper.
%\medskip
%\cred{
 %{\sc Note added at the revision}
%At the moment of finishing this revision, we became aware of the latest development on computing Tukey halfspace medain by
%M. Bogi\'cevi\'c and M. Merkle (2016).}

\end{document}